\documentstyle[12pt]{article}
\title{The algebraic de Rham theorem \\ for toric varieties}
\author{
{\sc Tadao Oda}\thanks{Partly supported by the Grants-in-Aid for
Co-operative Research as well as Scientific Research,
the Ministry of Education, Science and Culture, Japan.
\newline
\hspace*{1.5em} {\em $1991$ Mathematics Subject Classification\/}.
Primary 14M25;
Secondary 14F40, %de Rham cohomology
14F32,           %intersection cohomology
32S60.      %stratification, constructible sheaves, intersection cohomology
}}
\date{}
\newtheorem{Theorem}{{\sc Theorem}}[section]
\newtheorem{Proposition}[Theorem]{{\sc Proposition}}
\newtheorem{Lemma}[Theorem]{{\sc Lemma}}
\newtheorem{Corollary}[Theorem]{{\sc Corollary}}
\newcommand{\Proof}{{\sc Proof.} \quad}

\newcommand{\Remark}{{\sc Remark.} \quad}

\newcommand{\qed}{\hspace*{\fill}q.e.d.}
\newcommand{\C}{{\bf C}}
\newcommand{\Z}{{\bf Z}}
\newcommand{\R}{{\bf R}}
\newcommand{\Q}{{\bf Q}}
\newcommand{\bH}{{\bf H}}
\newcommand{\bP}{{\bf P}}

\newcommand{\cO}{{\cal O}}

\newcommand{\cX}{{\cal X}}
\newcommand{\cL}{{\cal L}}
\newcommand{\cD}{{\cal D}}
\newcommand{\cU}{{\cal U}}
\newcommand{\tDelta}{\widetilde{\Delta}}
\newcommand{\tPhi}{\widetilde{\Phi}}
\newcommand{\tX}{\widetilde{X}}
\newcommand{\tY}{\widetilde{Y}}
\newcommand{\oSigma}{\bar{\Sigma}}
\newcommand{\osigma}{\bar{\sigma}}
\newcommand{\otau}{\bar{\tau}}
\newcommand{\oN}{\bar{N}}
\newcommand{\on}{\bar{n}}
\newcommand{\oM}{\bar{M}}
\newcommand{\oeta}{\bar{\eta}}

\newcommand{\Hom}{\mathop{\rm Hom}\nolimits}

\newcommand{\Spec}{\mathop{\rm Spec}\nolimits}

\newcommand{\emb}{\mathop{\rm emb}\nolimits}
\newcommand{\orb}{\mathop{\rm orb}\nolimits}
\newcommand{\Lie}{\mathop{\rm Lie}\nolimits}

\newcommand{\Star}{\mathop{\rm Star}\nolimits}
\newcommand{\tcdot}{{\textstyle \cdot}}
\newcommand{\I}{\mbox{\rm {\scriptsize I}}}
\newcommand{\II}{\mbox{\rm {\scriptsize II}}}
\newcommand{\id}{{\rm id}}
\newcommand{\an}{^{{\rm an}}}
\newcommand{\tOmega}{\widetilde{\Omega}}
\begin{document}
\maketitle

\begin{center}
{\small Dedicated to Professor Takeshi Kotake on his sixtieth birthday}
\end{center}

%%%%%%%%%%%%%%%%%%%%%%%%%%%%%%%%%%%%%%%%%%%%%%%%%%%%%%%%%%%%%%%%%%%
\begin{abstract}
On an arbitrary toric variety, we introduce the logarithmic double complex,
which is essentially the same as the algebraic de Rham complex in the
nonsingular case, but which behaves much better in the singular case.

Over the field of complex numbers, we prove the toric analog of the
algebraic de Rham theorem which Grothendieck formulated and proved for
general nonsingular algebraic varieties re-interpreting an earlier work of
Hodge-Atiyah.
Namely, for a finite simplicial fan which need not be complete,
the complex cohomology groups of the corresponding
toric variety as an analytic space coincide with the hypercohomology
groups of the single complex associated to the logarithmic double complex.
They can then be described combinatorially as Ishida's cohomology groups
for the fan.

We also prove vanishing theorems for Ishida's cohomology groups.
As a consequence, we deduce directly that the complex cohomology groups
vanish in odd degrees for toric varieties which correspond to finite
simplicial fans with full-dimensional
convex support. In the particular case of
complete simplicial fans, we thus have a direct proof for an earlier
result of Danilov and the author.
\end{abstract}

%%%%%%%%%%%%%%%%%%%%%%%%%%%%%%%%%%%%%%%
\section*{Introduction}

Let $\Delta$ be a finite fan for a free $\Z$-module $N$ of rank $r$,
and denote by $X:=T_N\emb(\Delta)$ the associated $r$-dimensional
toric variety over the field $\C$ of complex numbers.
We also denote by $\cX:=X\an$ the associated complex analytic space.

In Section~1, we introduce the {\em logarithmic double complex\/}
$\cL_X^{\tcdot,\tcdot}$ of $\cO_X$-modules. The associated single complex
$\cL_X^{\tcdot}$, called the {\em logarithmic complex\/} for $X$, turns
out to be a natural toric generalization of the algebraic de Rham complex
for general nonsingular algebraic varieties.
Indeed, when $X$ is a nonsingular toric variety, then $\cL_X^{\tcdot}$
is canonically quasi-isomorphic to the algebraic de Rham complex
$\Omega_X^{\tcdot}$.
More generally, suppose $\Delta$ is a simplicial fan, so that $X$
has at worst quotient singularities. Then $\cL_X^{\tcdot}$ is
canonically quasi-isomorphic to the complex $\tOmega_X^{\tcdot}$ of
$\cO_X$-modules consisting of the Zariski differential forms on $X$.
By Danilov's Poincar\'e lemma \cite{danilov}
(see also Ishida \cite[Prop.~2.1]{ishida_derham}), we thus see that the
associated complex analytic complex $(\cL_X^{\tcdot})\an$ is
quasi-isomorphic to the constant sheaf $\C_{\cX}$.
We would like to point out here that
according to Ishida, $(\cL_X^{\tcdot})\an[2r]$ is quasi-isomorphic to
the {\em globally normalized dualizing complex\/} of $\C_{\cX}$-modules
in the sense of Verdier \cite{verdier}.

In Section 2, we recall the definition of Ishida's $p$-th complex
$C^{\tcdot}(\Delta,\Lambda^p)$ of $\Z$-modules for $0\leq p\leq r$
and its cohomology groups $H^q(\Delta,\Lambda^p)$,
which were considered in Ishida \cite{ishida_dualizing} and
\cite{odabook} in different notation.

Section 3 is devoted to the algebraic de Rham theorem for the logarithmic
complex $\cL_X^{\tcdot}$.
Without any condition on the fan $\Delta$, we first show that the
algebraic hypercohomology group $\bH^{\tcdot}(X,\cL_X^{\tcdot})$ is
a direct sum of the scalar extensions $H^{\tcdot}(\Delta,\Lambda^p)_{\C}$
to $\C$ of Ishida's cohomology groups for various $p$'s.

When $\Delta$ is simplicial but not necessarily complete, we mimic the
proof, due to Grothendieck \cite{grothendieck},
%and Hodge-Atiyah \cite{hodge_atiyah},
of the usual algebraic de Rham theorem
to show that the algebraic hypercohomology group
$\bH^{\tcdot}(X,\cL_X^{\tcdot})$ is canonically isomorphic to the
complex analytic hypercohomology group
$\bH^{\tcdot}(\cX,(\cL_X^{\tcdot})\an)$.
Consequently, we have our main result
\[
H^l(\cX,\C)=\bigoplus_{p+q=l}H^q(\Delta,\Lambda^p)_{\C}
\qquad\mbox{for each}\quad l
\]
for $\Delta$ simplicial but not necessarily complete.

In Section 4, we prove various vanishing theorems for Ishida's
cohomology groups.
When $\Delta$ is simplicial and complete, we have a direct combinatorial
proof for
\[
H^q(\Delta,\Lambda^p)_{\Q}=0\qquad\mbox{for}\quad q\neq p.
\]
Hence $H^l(\cX,\C)=0$ for $l$ odd, while
\[
H^{2p}(\cX,\C)=H^p(\Delta,\Lambda^p)_{\C}\qquad\mbox{for}\quad
0\leq p\leq r.
\]
In the complete nonsingular case, this was proved earlier by means of the
usual algebraic de Rham theorem by Danilov \cite{danilov}
and \cite[Thm.\ 3.11]{odabook}.

When $\Delta$ is simplicial with support $|\Delta|$ convex of dimension $r$
but not necessarily complete, we continue to have the same vanishing
theorem
\[
H^q(\Delta,\Lambda^p)_{\Q}=0\qquad\mbox{for}\quad q\neq p,
\]
which turns out to be a special case of a more general vanishing theorem
due to Ishida. In view of our algebraic de Rham theorem in Section~3,
we continue to have $H^l(\cX,\C)=0$ for $l$ odd, while
\[
H^{2p}(\cX,\C)=H^p(\Delta,\Lambda^p)_{\C}\qquad\mbox{for}\quad
0\leq p\leq r
\]
when $\Delta$ is simplicial with $|\Delta|$ convex of dimension $r$
but not necessarily complete.

The intersection cohomology theory of Goresky and MacPherson
\cite{gmI}, \cite{gmII}, which was further developed by
Beilinson, Bernstein and Deligne \cite{bbd}, attaches to
singular spaces new topological invariants much better behaved than
ordinary cohomology or homology groups.
Since a toric variety is a simple-minded singular space described
in terms of a fan, it is natural to try to describe its intersection
cohomology group and, more generally its intersection complex,
in terms of the fan as well.
We believe that the logarithmic complex
will play a key role in describing
the intersection cohomology group
(especially with respect to the middle perversity) of toric varieties.

There have been earlier attemps in this direction by
J. N. Bernstein, A. G. Khovanskii and R. D. MacPherson
(see Stanley \cite{stanley_aspm}), Denef and Loeser \cite{denef_loeser},
Fieseler \cite{fieseler} and others. However, they depend either on
the decomposition theorem of Beilinson, Bernstein, Deligne and Gabber
(see Beilinson, Bernstein and Deligne \cite{bbd}
as well as Goresky and MacPherson \cite{gm_map})
or on the purity theorem of Deligne and Gabber
(cf.\ Deligne \cite{deligne_weilII}, \cite{deligne_gabber}).

This problem turns out to be closely related to the problem of finding an
elementary proof of the strong Lefschetz theorem for toric varieties,
and hence of the elementary proof due to Stanley \cite{stanley} of the
``$g$-theorem'' for simplicial convex polytopes conjectured earlier by
McMullen \cite{mcmullen}.
(See also \cite{lefsch}, \cite{hydera},
Stanley \cite{stanley_local_h} and
McMullen \cite{mcmullen2}.)

We refer the reader to \cite{odabook}
for basic results on toric varieties used in this paper.

Thanks are due to Dr.~M.-N.~Ishida for stimulating discussion.

%%%%%%%%%%%%%%%%%%%%%%%%%%%%%%%%%%%%%%%%%%%%%%%%%%%

\section{Logarithmic double complex}

Let $\Delta$ be a finite fan for a free $\Z$-module $N$ of rank $r$,
and denote by $X:=T_N\emb(\Delta)$ the associated $r$-dimensional
toric variety over the field $\C$ of complex numbers.
The complement $D:=X\setminus T_N$ of the algebraic torus
$T_N:=N\otimes_{\Z}\C^*\cong(\C^*)^r$ is a Weil divisor on $X$ but is
not necessarily a Cartier divisor.
The dual $\Z$-module $M:=\Hom_{\Z}(N,\Z)$, with the canonical bilinear
pairing $\langle\phantom{m},\phantom{m}\rangle :M\times N\rightarrow\Z$,
is isomorphic to the character group of the algebraic torus $T_N$.
For each $m\in M$ we denote the corresponding character by
$t^m:T_N\rightarrow\C^*$ (which was denoted by ${\bf e}(m)$ in \cite{odabook}),
and identify $\C[M]:=\oplus_{m\in M}\C t^m$ with the group algebra of $M$
over $\C$. Hence $T_N$ is the group of $\C$-valued points of the
group scheme $\Spec(\C[M])$.

Each $n\in N$ gives rise to a $\C$-derivation $\delta_n$ of $\C[M]$
defined by $\delta_n(t^m):=\langle m,n\rangle t^m$. Consequently, we
have a canonical isomorphism to the Lie algebra
\[
\C\otimes_{\Z}N\stackrel{\sim}{\longrightarrow}\Lie(T_N),
\qquad 1\otimes n\mapsto\delta_n,
\]
hence an $\cO_X$-isomorphism
$\cO_X\otimes_{\Z}N\stackrel{\sim}{\rightarrow}\Theta_X(-\log D)$, where the
right hand side is the sheaf of germs of algebraic vector fields on $X$
with logarithmic zeros along the Weil divisor $D$. Its dual
$\Omega^1_X(\log D)$ is the sheaf of germs of algebraic $1$-forms with
logarithmic poles along $D$, and we get an $\cO_X$-isomorphism
\[
\cO_X\otimes_{\Z}M\stackrel{\sim}{\longrightarrow}\Omega^1_X(\log D),
\qquad 1\otimes m\mapsto \frac{dt^m}{t^m}.
\]
Taking the exterior product, we thus get an $\cO_X$-isomorphism
$\cO_X\otimes_{\Z}\bigwedge^{\tcdot}M\stackrel{\sim}{\rightarrow}
\Omega^{\tcdot}_X(\log D)$. The exterior differentiation $d$ on the
right hand side corresponds to the operation on the left hand side
which sends a $p$-form $t^m\otimes m_1\wedge\cdots\wedge m_p$ to the
$(p+1)$-form $t^m\otimes m\wedge m_1\wedge\cdots\wedge m_p$
(cf. \cite[Chap.~3]{odabook}).

Recall that the set of $T_N$-orbits in $X$ is in one-to-one correspondence
with $\Delta$ by the map which sends each $\sigma\in\Delta$ to the
$T_N$-orbit
\[
\orb(\sigma)=\Spec(\C[M\cap\sigma^{\perp}])=T_{N/\Z(N\cap\sigma)}.
\]
The closure $V(\sigma)$ in $X$ of $\orb(\sigma)$ is known to be a toric
variety with respect to a fan for the $\Z$-module $N/\Z(N\cap\sigma)$.
Namely,
\[
V(\sigma)=
T_{N/\Z(N\cap\sigma)}\emb(\{(\tau+(-\sigma))/\R\sigma\;\mid\;
\tau\in\Star_{\sigma}(\Delta)\}),
\]
where $\R\sigma=\sigma+(-\sigma)$ is the smallest $\R$-subspace containing
$\sigma$ of $N_{\R}:=N\otimes_{\Z}\R$, while
$\Star_{\sigma}(\Delta):=\{\tau\in\Delta\mid \tau\succ\sigma\}$.
Hence $D(\sigma):=V(\sigma)\setminus\orb(\sigma)$ is a Weil divisor
on $V(\sigma)$. In particular, we have
$\orb(\{0\})=T_N$, $D(\{0\})=D$ and $V(\{0\})=X$.

For each integer $q$ with $0\leq q\leq r$, denote
$\Delta(q):=\{\sigma\in\Delta\mid \dim\sigma=q\}$. For each pair of
integers $p$, $q$, let
\[
\cL^{p,q}_X:=\bigoplus_{\sigma\in\Delta(q)}
\Omega^{p-q}_{V(\sigma)}(\log D(\sigma))
=\bigoplus_{\sigma\in\Delta(q)}\cO_{V(\sigma)}\otimes_{\Z}
\bigwedge\nolimits^{p-q}(M\cap\sigma^{\perp})\qquad\mbox{if}\quad
0\leq q\leq p,
\]
and $\cL^{p,q}=0$ otherwise.
$d_{\I}:\cL^{p,q}_X\rightarrow\cL^{p+1,q}_X$ is defined to be the direct
sum of the exterior differentiation for each $\sigma\in\Delta(q)$.
We define $d_{\II}:\cL^{p,q}_X\rightarrow\cL^{p,q+1}_X$ as follows:
The $(\sigma,\tau)$-component of $d_{\II}$ for $\sigma\in\Delta(q)$ and
$\tau\in\Delta(q+1)$ is defined to be zero when $\sigma$ is not a face of
$\tau$. On the other hand, if $\sigma$ is a face of $\tau$, then
a primitive element $n\in N$ is uniquely determined modulo $N\cap\R\sigma$
so that $\tau+(-\sigma)=\R_{\geq 0}n+\R\sigma$.
$M\cap\tau^{\perp}$ is a $\Z$-submodule of corank one
in $M\cap\sigma^{\perp}$. The $(\sigma,\tau)$-component of $d_{\II}$ in
this case is then defined to be the tensor product of the restriction
homomorphism $\cO_{V(\sigma)}\rightarrow\cO_{V(\tau)}$ with the
interior product with respect to $n$. Namely, an element in the
$\sigma$-component of $\cL^{p,q}_X$ of the form
\[
t^m\otimes m_1\wedge m_2\wedge\cdots\wedge m_{p-q},\qquad
m, m_1\in M\cap\sigma^{\perp},\quad
m_2,\ldots,m_{p-q}\in M\cap\tau^{\perp}
\]
is sent to $t^m\otimes\langle m_1,n\rangle m_2\wedge\cdots\wedge m_{p-q}$
if $m\in M\cap\tau^{\perp}$, and to 0 otherwise.
$d_{\II}$ is the Poincar\'e residue map.

$d_{\I}\circ d_{\I}=0$ and $d_{\I}\circ d_{\II}+d_{\II}\circ d_{\I}=0$
are obvious, while $d_{\II}\circ d_{\II}=0$
was shown in Ishida \cite[Lemma 1.4 and Prop. 1.6]{ishida_dualizing}.
Consequently, we get a double complex $\cL^{\tcdot,\tcdot}_X$
of $\cO_X$-modules,
which we call the {\em logarithmic double complex\/} for the toric variety
$X$. The associated single complex is denoted by $\cL^{\tcdot}_X$ and
is called the {\em logarithmic complex\/} for $X$.

For simplicity, we denote by $\cX:=X\an$ the complex analytic space
associated to the toric variety $X=T_N\emb(\Delta)$.

\begin{Theorem} \label{thm_simplicialLdotan}
If $\Delta$ is a {\em simplicial\/} fan, then we have a quasi-isomorphism
$\C_{\cX}\simeq(\cL^{\tcdot}_X)\an$.
\end{Theorem}

\Proof
It suffices to prove the assertion when $X$ is affine, and there is an
easy proof in that case. However, we here give a direct proof valid for
general $X$.

Let $\tOmega^{\tcdot}_X:=j_*\Omega^{\tcdot}_U$
be the complex of $\cO_X$-modules consisting of the Zariski differential
forms on $X$, where $j:U\rightarrow X$ is the open immersion of the
smooth locus $U$ of $X$. By Danilov's Poincar\'e lemma
\cite{danilov} (see also Ishida \cite[Prop.~2.1]{ishida_derham}),
we have a quasi-isomorphism $\C_{\cX}\simeq(\tOmega^{\tcdot}_X)\an$
for an {\em arbitrary\/} fan $\Delta$. Furthermore, if $\Delta$ is
simplicial, then by \cite[Theorem 3.6]{odabook} we have a quasi-isomorphism
$\tOmega^p_X\simeq\cL^{p,\tcdot}_X$ for each fixed $p$.
\qed
\bigskip

\Remark
In \cite{odabook}, $\cL^{p,\tcdot}_X$ was denoted by ${\cal K}^{\tcdot}(X;p)$.

As for an {\em arbitrary\/} fan $\Delta$ for $N\cong(\Z)^r$ which need not be
simplicial, Ishida has a proof for the following amazing result:
$(\cL^{\tcdot}_X)\an[2r]$ is quasi-isomorphic to the {\em globally normalized
dualizing complex\/} $\cD^{\tcdot}_{\cX}$ of $\C_{\cX}$-modules
in the sense of Verdier \cite{verdier}.
The point of this result of Ishida's lies in the fact that the dualizing
complex $\cD^{\tcdot}_{\cX}$ can be expressed in terms of a complex
comprizing of {\em algebraic\/} and {\em coherent\/} $\cO_X$-modules.
Analogously, Ishida \cite[Theorem 3.3]{ishida_dualizing} and
\cite[Theorem 5.4]{ishida_derham} earlier showed $\cL^{r,\tcdot}_X[r]$
to be quasi-isomorphic to the {\em globally normalized dualizing complex\/}
of $\cO_X$-modules.
\bigskip

%%%%%%%%%%%%%%%%%%%%%%%%%%%%%%%%%%%%%%%%%%%%%%%%%%
The following is a toric generalization of a result due to Grothendieck
and Deligne \cite[I, \S 3 and \S 6]{deligne_equadiff} on the complement
in a smooth variety of a divisor with normal crossings.

\begin{Proposition} \label{prop_Rj*CTN}
Let $j$ be the open immersion of $T_N$ into a toric variety $X$, and
denote by $j\an:(T_N)\an\rightarrow X\an$ the corresponding open
immersion of complex analytic spaces. Then we have a canonical
quasi-isomorphism
\[
\R j\an_*\C_{(T_N)\an}\simeq\left(\Omega^{\tcdot}_X(\log D)\right)\an
=(\cO_{X})\an\otimes_{\Z}\bigwedge\nolimits^{\tcdot}M
\]
and canonical isomorphisms
\[
H^{\tcdot}((T_N)\an,\C)=\bH^{\tcdot}(X,\Omega^{\tcdot}_X(\log D))
=\C\otimes_{\Z}\bigwedge\nolimits^{\tcdot}M,
\]
where $\bH^{\tcdot}(X,\Omega^{\tcdot}_X(\log D))$ is the hypercohomology
group of the complex $\Omega^{\tcdot}_X(\log D)$ consisting of $\cO_X$-modules.
\end{Proposition}

\Proof
We may assume $X$ to be smooth, since
$\Omega^{\tcdot}_X(\log D)=\cO_X\otimes_{\Z}\bigwedge^{\tcdot}M$
and since toric singularities are rational so that
$\R f_*\cO_{X'}=\cO_X$ holds for any equivariant resolution
of singularities $f:X'\rightarrow X$
(see, for instance, \cite[Cor.\ 3.9]{odabook} and
Ishida \cite[Cor.\ 3.3]{ishida_derham}).
Consequently, $D:=X\setminus T_N$ is a divisor with simple normal
crossings on a smooth $X$ as in the situation dealt with by
Grothendieck and Deligne. Let us repeat their proof here for convenience.

For the proof of the first assertion, we may obviously assume $X$ to be affine
as well, hence $j\an$ is a Stein morphism. By the Poincar\'e lemma,
we thus have
\[
\R j_*\an\C_{(T_N)\an}=\R j_*\an(\Omega^{\tcdot}_{T_N})\an
=j_*\an(\Omega^{\tcdot}_{T_N})\an.
\]
The term on the extreme right hand side is canonically quasi-isomorphic
to $(\Omega^{\tcdot}_X(\log D))\an$ by
Deligne \cite[Lemma 6.9]{deligne_equadiff}.

As for the second assertion, note that the equality
$H^{\tcdot}((T_N)\an,\C)=\bigwedge^{\tcdot}M_{\C}$ itself is well-known,
where $M_{\C}:=M\otimes_{\Z}\C$.
To show the whole set of the equalities, we start with the consequence
\[
H^{\tcdot}((T_N)\an,\C)=\bH^{\tcdot}(X\an,\R j_*\an\C_{(T_N)\an})
=\bH^{\tcdot}(X\an,(\Omega_X^{\tcdot}(\log D))\an)
\]
of the first assertion. The term on the extreme right hand side is canonically
isomorphic to the algebraic hypercohomology group
$\bH^{\tcdot}(X,\Omega^{\tcdot}_X(\log D))$ by
Deligne \cite[Thm.\ 6.2]{deligne_equadiff}.
Since $\Omega_X^{\tcdot}(\log D)=\cO_X\otimes_{\Z}\bigwedge^{\tcdot}M$
has $\cO_X$-coherent components,
the algebraic hypercohomology group in question coincides with the cohomology
group of the single complex associated to the \v{C}ech double complex
$\check{C}^{\tcdot}(\cU,\Omega_X^{\tcdot}(\log D))$ with respect
to the $T_N$-stable
affine open covering $\cU:=\{U_{\sigma}\mid\sigma\in\Delta\}$
with $U_{\sigma}:=\Spec(\C[M\cap\sigma^{\vee}])$. Moreover,
$T_N$ has a canonical algebraic action on the \v{C}ech double complex,
which gives rise to an eigenspace decomposition
\[
\check{C}^{\tcdot}(\cU,\Omega_X^{\tcdot}(\log D))=
\bigoplus_{m\in M}\check{C}^{\tcdot}(\cU,\Omega_X^{\tcdot}(\log D))_m
\]
with respect to the characters $m\in M$ of $T_N$.
As a result, we have an eigenspace decomposition
\[
\bH^{\tcdot}(X,\Omega_X^{\tcdot}(\log D))
=\bigoplus_{m\in M}\bH^{\tcdot}(X,\Omega_X^{\tcdot}(\log D))_m
\]
for the hypercohomology group as well.

For $m\neq 0$ we have $\bH^{\tcdot}(X,\Omega_X^{\tcdot}(\log D))_m=0$,
since the $m$-th component $d_m$ of the exterior differentiation
is the exterior multiplication by $m$, hence is exact.

On the other hand, for $m=0$ we have $d_0=0$, hence the \v{C}ech
double complex is a direct sum of the complexes
$\check{C}^{\tcdot}(\cU,\Omega_X^p(\log D))_0$ for $0\leq p\leq r$.
Since $\Omega_X^p(\log D)=\cO_X\otimes_{\Z}\bigwedge^p M$ and since
\[
H^l(X,\cO_X)_0=
\left\{\begin{array}{lll}
\C & & l=0     \\
0  & & l\neq 0. \\
\end{array}\right.
\]
as we recall later in Lemma \ref{lem_HVOV0}, we conclude that
\[
\bH^{\tcdot}(X,\Omega_X^{\tcdot}(\log D))_0
=H^0(X,\Omega_X^{\tcdot}(\log D))_0=\bigwedge\nolimits^{\tcdot}M_{\C}.
\]
\qed
\bigskip

Applying Proposition~\ref{prop_Rj*CTN} to the immersion
$j_{\sigma}:\orb(\sigma)\rightarrow X$ for each $\sigma\in\Delta$, we
get the following:

\begin{Corollary} \label{cor_Ldotqan}
For each fixed $q$, we have a canonical quasi-isomorphism
\[
\left(\cL^{\tcdot,q}_X\right)\an\simeq
\bigoplus_{\sigma\in\Delta(q)}\R(j_{\sigma})\an_*\C_{\orb(\sigma)}[-q],
\]
where $[-q]$ denotes the degree shift to the right by $q$.
\end{Corollary}

%%%%%%%%%%%%%%%%%%%%%%%%%%%%%%%%%%%%%%%%%%%%%%%%%%%

\section{Ishida's complexes}

Let $\Delta$ be a finite fan for a free $\Z$-module $N$ of rank $r$,
and for $0\leq q\leq r$ denote
$\Delta(q):=\{\sigma\in\Delta\mid\dim\sigma=q\}$ as before.

For each integer $p$ with $0\leq p\leq r$, {\em Ishida's $p$-th complex\/}
$C^{\tcdot}(\Delta,\Lambda^p)$ of $\Z$-modules is defined as follows:
\[
C^q(\Delta,\Lambda^p):=\bigoplus_{\sigma\in\Delta(q)}
\bigwedge\nolimits^{p-q}(M\cap\sigma^{\perp})\qquad\mbox{if}\quad
0\leq q\leq p,
\]
and $C^q(\Delta,\Lambda^p)=0$ otherwise.
For $\sigma\in\Delta(q)$ and $\tau\in\Delta(q+1)$,
the $(\sigma,\tau)$-component of the coboundary map
\[
\delta:C^q(\Delta,\Lambda^p)=\bigoplus_{\sigma\in\Delta(q)}
\bigwedge\nolimits^{p-q}(M\cap\sigma^{\perp})\longrightarrow
C^{q+1}(\Delta,\Lambda^p)=\bigoplus_{\tau\in\Delta(q+1)}
\bigwedge\nolimits^{p-q-1}(M\cap\tau^{\perp})
\]
is defined to be $0$ if $\sigma$ is not a face of $\tau$.
On the other hand, if $\sigma$ is a face of $\tau$, then a primitive
element $n\in N$ is uniquely determined modulo $N\cap\R\sigma$ so that
$\tau+(-\sigma)=\R_{\geq 0}n+\R\sigma$. The $(\sigma,\tau)$-component
of $\delta$ in this case is defined to be the interior product with
respect to this $n$. Namely, the element
$m_1\wedge m_2\wedge\cdots\wedge m_{p-q}$ with $m_1\in M\cap\sigma^{\perp}$
and $m_2,\ldots,m_{p-q}\in M\cap\tau^{\perp}$ is sent to
$\langle m_1,n\rangle m_2\wedge\cdots\wedge m_{p-q}$.
As Ishida \cite[Prop.\ 1.6]{ishida_dualizing} showed, $\delta\circ\delta=0$
holds so that $C^{\tcdot}(\Delta,\Lambda^p)$ is a complex of $\Z$-modules.
We denote its cohomology group by $H^{\tcdot}(\Delta,\Lambda^p)$.
We will be mainly concerned with their scalar extensions
$C^{\tcdot}(\Delta,\Lambda^p)_{\Q}$, $C^{\tcdot}(\Delta,\Lambda^p)_{\C}$,
$H^{\tcdot}(\Delta,\Lambda^p)_{\Q}$, $H^{\tcdot}(\Delta,\Lambda^p)_{\C}$
to $\Q$ and $\C$.

By definition, we have $H^q(\Delta,\Lambda^p)=0$ unless $0\leq q\leq p$.
\bigskip

\Remark
In \cite[\S 3.2]{odabook}, $C^{\tcdot}(\Delta,\Lambda^p)$
and $H^{\tcdot}(\Delta,\Lambda^p)$ were denoted by $C^{\tcdot}(\Delta;p)$ and
$H^{\tcdot}(\Delta;p)$, respectively.

As in \cite{lefsch}, we can define a similar complex
$C^{\tcdot}(\Pi,{\cal G}_p)$ of $\R$-vector spaces for a {\em simplicial\/}
polyhedral cone decomposition $\Pi$ of an $\R$-vector space
endowed with a marking for each one-dimensional cone in $\Pi$.
Note, however, that unless a lattice $N$ is given as in the case of a fan,
we cannot define the coboundary map in the case of a {\em non-simplicial\/}
convex polyhedral cone decomposition $\Pi$ even if it is endowed with a
marking. It is crucial that for a codimension one face $\sigma$ of
$\tau$ in a fan, a primitive element $n\in N$ is uniquely determined
modulo $N\cap\R\sigma$ so that $\tau+(-\sigma)=\R_{\geq 0}n+\R\sigma$
holds as above, regardless of whether $\tau$ is simplicial or not.
\bigskip

The following result is slightly stronger than \cite[Lemma 3.7]{odabook},
and the proof is similar to that for \cite[Prop.\ 3.5]{lefsch},
which concerns analogous $\R$-coefficient cohomology groups for
the simplicial polyhedral cone decomposition consisting of the
faces of a simplicial cone in a finite dimensional $\R$-vector space:

\begin{Proposition} \label{prop_vanishing_simplicialcone}
Let $\pi$ be a {\em simplicial} rational polyhedral cone in $N_{\R}$.
Then for each $0\leq p\leq r$, the cohomology group of Ishida's $p$-th
complex for the
fan $\Gamma_{\pi}$ consisting of all the faces of $\pi$ satisfies
\[
H^q(\Gamma_{\pi},\Lambda^p)_{\Q}:=
H^q(\Gamma_{\pi},\Lambda^p)\otimes_{\Z}\Q=
\left\{
\begin{array}{lll}
\bigwedge^p(M_{\Q}\cap\pi^{\perp}) &\phantom{mm}& q=0 \\
0                                  &            & q\neq 0, \\
\end{array}
\right.
\]
where $M_{\Q}:=M\otimes_{\Z}\Q$.
\end{Proposition}

%%%%%%%%%%%%%%%%%%%%%%%%%%%%%%%%%%%%%%%%%%%%%%%%%%%

\section{The algebraic de Rham theorem}

In this section, we denote by $X:=T_N\emb(\Delta)$ the $r$-dimensional
toric variety over $\C$ corresponding to a finite fan $\Delta$ for
$N\cong\Z^r$.
For simplicity, we again denote the corresponding complex analytic space
by $\cX:=X\an$.

\begin{Proposition} \label{prop_alglog_ishida}
For an {\em arbitrary} fan $\Delta$ which need not be complete nor simplicial,
the hypercohomology group of the logarithmic complex
$\cL^{\tcdot}_X$ has a direct sum decomposition
\[
\bH^l(X,\cL^{\tcdot}_X)=\bigoplus_{p+q=l}H^q(\Delta,\Lambda^p)_{\C}
\qquad\mbox{for each $l$}.
\]
\end{Proposition}

\Proof
Consider the $T_N$-stable
affine open covering $\cU:=\{U_{\sigma}\mid\sigma\in\Delta\}$
of $X$ with $U_{\sigma}:=\Spec(\C[M\cap\sigma^{\vee}])$.
We know that (cf.\ \cite{odabook})
\[
U_{\sigma_0}\cap U_{\sigma_1}\cap\cdots\cap U_{\sigma_q}
=U_{\sigma_0\cap\sigma_1\cap\cdots\cap\sigma_q}\qquad\mbox{for all}\quad
\sigma_0,\sigma_1,\ldots,\sigma_q\in\Delta.
\]
Since each component of $\cL^{\tcdot}_X$ is a coherent $\cO_X$-module,
the hypercohomology group $\bH^{\tcdot}(X,\cL^{\tcdot}_X)$ coincides with
the cohomology group of the single complex associated to the
\v{C}ech double complex $\check{C}^{\tcdot}(\cU,\cL^{\tcdot}_X)$ with
respect to the affine open covering $\cU$.

Because of the canonical algebraic action of the algebraic torus $T_N$ on
$\check{C}^{\tcdot}(\cU,\cL^{\tcdot}_X)$, we have the eigenspace
decomposition
\[
\check{C}^{\tcdot}(\cU,\cL_X^{\tcdot})
=\bigoplus_{m\in M}\check{C}^{\tcdot}(\cU,\cL_X^{\tcdot})_m
\qquad\mbox{hence}\qquad
\bH^{\tcdot}(\cU,\cL_X^{\tcdot})
=\bigoplus_{m\in M}\bH^{\tcdot}(\cU,\cL_X^{\tcdot})_m
\]
with respect to the characters $m\in M$.
Let us consider the triple complex
\[
\check{C}_m^{\tcdot,\tcdot,\tcdot}:=
\check{C}^{\tcdot}(\cU,\cL_X^{\tcdot,\tcdot})_m\qquad\mbox{with}\quad
\check{C}_m^{l,p,q}:=
\check{C}^{l}(\cU,\cL_X^{p,q})_m.
\]
The differential for $l$ is the $m$-th component $\check{\delta}_m$
of the \v{C}ech coboundary $\check{\delta}$, while those for $p$ and
$q$ are the $m$-th components $(d_{\I})_m$ and $(d_{\II})_m$ of
$d_{\I}$ and $d_{\II}$, respectively.

$\bH^{\tcdot}(X,\cL_X^{\tcdot})_m$ is the cohomology group of the
single complex associated to $\check{C}^{\tcdot,\tcdot,\tcdot}_m$.
For $m\neq 0$, we have $\bH^{\tcdot}(X,\cL_X^{\tcdot})_m=0$, hence
$\bH^{\tcdot}(X,\cL_X^{\tcdot})=\bH^{\tcdot}(X,\cL_X^{\tcdot})_0$
is $T_N$-invariant.
Indeed, since $(d_{\I})_m$ is the exterior multiplication by $m$,
the complex $(\check{C}_m^{l,\tcdot,q},(d_{\I})_m)$ is acyclic for
all $l$ and $q$.
Consequently, the cohomology group vanishes for
the single complex associated to the triple complex
$\check{C}^{\tcdot,\tcdot,\tcdot}_m$.

On the other hand, for $m=0$ we have $(d_{\I})_0=0$, hence
$\check{C}^{\tcdot,\tcdot,\tcdot}_0=
\bigoplus_p\check{C}^{\tcdot,p,\tcdot}_0$ is a direct sum of double
complexes.
Since $\cL_X^{p,q}$ is a direct sum of coherent sheaves of the form
$\cO_{V(\sigma)}$, we see by Lemma \ref{lem_HVOV0} below that
the cohomology group of the single complex associated to
$\check{C}^{\tcdot,p,\tcdot}_0$ coincides with that of
$H^0(X,\cL_X^{p,\tcdot})_0=C^{\tcdot}(\Delta,\Lambda^p)_{\C}$.
Consequently, the cohomology group
$\bH^{\tcdot}(X,\cL_X^{\tcdot})=\bH^{\tcdot}(X,\cL_X^{\tcdot})_0$
of the single complex associated to $\check{C}_0^{\tcdot,\tcdot,\tcdot}$
is the direct sum of $H^{\tcdot}(\Delta,\Lambda^p)_{\C}$ as claimed.
\qed

\begin{Lemma} \label{lem_HVOV0}
For any toric variety $V$ with respect to an algebraic torus $T_N$,
the $T_N$-invariant part of the cohomology group for the structure sheaf
$\cO_V$ satisfies
\[
H^l(V,\cO_V)_0=
\left\{\begin{array}{lll}
\C & & l=0     \\
0  & & l\neq 0. \\
\end{array}\right.
\]
\end{Lemma}

\Proof
Demazure and Danilov gave a general description of the eigenspace
$H^{\tcdot}(V,L)_m$, with respect to a character $m$, of the $T_N$-action on
the cohomology group of an equivariant line bundle $L$ on $V$ in terms of the
corresponding support function (see, for instance, \cite[Thm.\ 2.6]{odabook}).
The present lemma is the special case $L=\cO_V$ and $m=0$.
\qed
\bigskip

We are now ready to state a generalization, in the toric context, of
the algebraic de Rham theorem due to Grothendieck \cite{grothendieck}.
An algebro-geometric proof valid in the case of complete nonsingular fans
can be found in Danilov \cite{danilov} and in \cite[Theorem 3.11]{odabook}.

\begin{Theorem} \label{thm_algderham}
{\rm (The algebraic de Rham theorem)}
For a {\em simplicial} fan $\Delta$ which need not be complete, we have the
following for each $l$:
\[
H^l(\cX,\C)=\bH^l(\cX,(\cL^{\tcdot}_X)\an)
=\bH^l(X,\cL^{\tcdot}_X)=\bigoplus_{p+q=l}H^q(\Delta,\Lambda^p)_{\C},
\]
where the second term from the left is the corresponding analytic
hypercohomology group.
\end{Theorem}

\Proof
The first equality follows from Theorem \ref{thm_simplicialLdotan},
while the third is a consequence of Proposition \ref{prop_alglog_ishida}.
To show the second equality, we mimic the proof due to
Grothendieck \cite{grothendieck}
in the case of the usual de Rham complex on smooth algebraic varieties.

There is a canonical homomorphism from the spectral sequence
\[
E_1^{p,q}:=H^q(X,\cL_X^p)\Longrightarrow
\bH^{p+q}(X,\cL_X^{\tcdot})
\]
for the algebraic hypercohomology to the spectral sequence
\[
E_1^{p,q}:=H^q(\cX,(\cL_X^p)\an)\Longrightarrow
\bH^{p+q}(\cX,(\cL_X^{\tcdot})\an)
\]
for the complex analytic hypercohomology.
To show that the homomorphism between the hypercohomology groups is an
isomorphism, we may assume $X$ to be affine, since the algebraic
(resp.\ complex analytic) hypercohomology group is the cohomology group
of the single complex associated to the \v{C}ech double complex
with respect to the affine open covering $\cU$
(resp.\ the corresponding Stein open covering $(\cU)\an$)
as in the proof of Proposition \ref{prop_alglog_ishida}.

We thus assume $X=U_{\pi}:=\Spec(\C[M\cap\pi^{\vee}])$ for a simplicial
rational cone $\pi\subset N_{\R}$, hence $\Delta=\Gamma_{\pi}$.

By Propositions \ref{prop_alglog_ishida} and
\ref{prop_vanishing_simplicialcone}, we have
\[
\bH^l(U_{\pi},\cL^{\tcdot}_{U_{\pi}})
=\bigoplus_{p+q=l}H^q(\Gamma_{\pi},\Lambda^p)_{\C}
=\bigwedge\nolimits^l(M\cap\pi^{\perp})_{\C}
\]
for each $l$.
We are done in view of
\[
\bH^l((U_{\pi})\an,(\cL^{\tcdot}_{U_{\pi}})\an)
=H^l((U_{\pi})\an,\C)=\bigwedge\nolimits^l(M\cap\pi^{\perp})_{\C}
\]
by Theorem \ref{thm_simplicialLdotan} and Danilov \cite[Lemma 12.3]{danilov}.
\qed
\bigskip

\Remark
Suppose $\Delta$ is an {\em arbitrary\/} finite fan which need not be
simplicial nor complete.
In view of the Verdier duality \cite{verdier} and Ishida's result mentioned
in the remark after Theorem \ref{thm_simplicialLdotan}, we have
\[
\bH^l(\cX,(\cL^{\tcdot}_X)\an)=\bH^l(\cX,\cD^{\tcdot}_{\cX}[-2r])
=\Hom_{\C}(H_c^{2r-l}(\cX,\C),\C),
\]
where the extreme right hand side is the dual of the cohomology group with
compact support.
Consequently, if we can show
\[
\Hom_{\C}(H_c^{2r-l}((U_{\pi})\an,\C),\C)
=\bigoplus_{p+q=l}H^q(\Gamma_{\pi},\Lambda^p)_{\C}
\qquad\mbox{with}\quad U_{\pi}:=T_N\emb(\Gamma_{\pi})
\]
for an arbitrary strongly convex rational polyhedral cone $\pi\subset N_{\R}$
which need not be simplicial, then Theorem \ref{thm_algderham}
could be generalized for an {\em arbitrary\/} fan in the form
\[
\Hom_{\C}(H_c^{2r-l}(\cX,\C),\C)=\bH^l(\cX,(\cL^{\tcdot}_X)\an)
=\bH^l(X,\cL^{\tcdot}_X)
=\bigoplus_{p+q=l}H^q(\Delta,\Lambda^p)_{\C}.
\]
\bigskip

%%%%%%%%%%%%%%%%%%%%%%%%%%%%%%%%%%%%%%%%%%%%%%%%%%%

\section{Vanishing theorems}

Analogs of the following were proved in Danilov \cite{danilov} and
\cite[Theorem 3.11]{odabook}
for complete nonsingular fans by means of the algebraic de Rham theorem,
and then directly in \cite[Theorem 4.1]{lefsch} for
complete simplicial polyhedral cone decompositions endowed with markings:

\begin{Proposition} \label{prop_vanishing_complete}
Let $\Delta$ be a {\em simplicial} and {\em complete} fan for
$N\cong\Z^r$. Then for all $p$ with $0\leq p\leq r$ we have
\[
H^q(\Delta,\Lambda^p)_{\Q}=0\qquad\mbox{for}\quad q\neq p.
\]
Moreover for all $p$, we have a perfect pairing
\[
H^p(\Delta,\Lambda^p)_{\Q}\times H^{r-p}(\Delta,\Lambda^{r-p})_{\Q}
\longrightarrow H^r(\Delta,\Lambda^r)_{\Q}\cong\bigwedge\nolimits^r M_{\Q}.
\]
Consequently, the corresponding toric variety $X:=T_N\emb(\Delta)$
satisfies $H^l(X\an,\C)=0$ for odd $l$, while
\[
H^{2p}(X\an,\C)=H^p(\Delta,\Lambda^p)_{\C}\qquad\mbox{for}\quad
0\leq p\leq r.
\]
\end{Proposition}

\Proof
We just indicate necessary modifications of the proof of
\cite[Theorem 4.1]{lefsch}.

We introduce a first quadrant double complex
$(K^{\tcdot,\tcdot},d',d'')$ as follows:
For nonnegative integers $i$ and $j$, let
\[
K^{i,j}:=\bigoplus_{\varphi\in\Delta(r-i)}\;\;
\bigoplus_{\stackrel{\scriptstyle \sigma\in\Delta(j)}{\sigma\prec\varphi}}
\bigwedge\nolimits^{p-j}
\left(M\cap\sigma^{\perp}\right)
\otimes_{\Z}(\det\varphi)^{-1},
\]
where $\det\varphi:=\bigwedge^{\dim\varphi}(N\cap\R\varphi)$
is the orientation $\Z$-module of rank one, and
$(\det\varphi)^{-1}$ is its dual $\Z$-module.
If $\psi\in\Delta(k-1)$ is a facet of $\varphi\in\Delta(k)$, we have
mutually dual nonzero orientation $\Z$-homomorphisms
\[
\det\psi\longrightarrow\det\varphi\qquad\mbox{and}\qquad
(\det\varphi)^{-1}\longrightarrow(\det\psi)^{-1}
\]
in the following manner: A primitive element $n\in N$ is uniquely
determined modulo $N\cap\R\psi$ so that
$\varphi+(-\psi)=\R_{\geq 0}n+\R\psi$. Then
$n_1\wedge\cdots\wedge n_{k-1}\in\det\psi$ is sent to
$n\wedge n_1\wedge\cdots\wedge n_{k-1}\in\det\varphi$.

For $\Delta(j)\ni\sigma\prec\varphi\in\Delta(r-i)$ and
$\Delta(j)\ni\sigma'\prec\psi\in\Delta(r-i-1)$, we define the
component of $d':K^{i,j}\rightarrow K^{i+1,j}$ with respect to
$(\varphi,\sigma)$ and $(\psi,\sigma')$ to be nonzero only when
$\varphi\succ\psi\succ\sigma=\sigma'$ and to be equal to
$(-1)^{j}$ times the dual orientation $\Z$-homomorphism
$(\det\varphi)^{-1}\rightarrow(\det\psi)^{-1}$ tensored with the
identity map for $\bigwedge^{p-j}(M\cap\sigma^{\perp})$.
On the other hand, for $\Delta(j)\ni\sigma\prec\varphi\in\Delta(r-i)$
and $\Delta(j+1)\ni\tau\prec\varphi'\in\Delta(r-i)$, we define the
component of $d'':K^{i,j}\rightarrow K^{i,j+1}$ with respect to
$(\varphi,\sigma)$ and $(\varphi',\tau)$ to be nonzero only when
$\varphi=\varphi'\succ\tau\succ\sigma$ and to be equal to the
tensor product of the homomorphism
$\bigwedge^{p-j}(M\cap\sigma^{\perp})\rightarrow
\bigwedge^{p-j-1}(M\cap\tau^{\perp})$ appearing in the definition of
Ishida's $p$-th complex, with the identity map for $(\det\varphi)^{-1}$.
It is easy to show that $(d')^2=(d'')^2=d'd''+d''d'=0$, hence
$K^{\tcdot,\tcdot}$ is a double complex of $\Z$-modules.

Its scalar extension $K^{\tcdot,\tcdot}_{\Q}$ to $\Q$ turns out to
satisfy
\[
H_{\I}^i(K^{\tcdot,j}_{\Q})=
\left\{\begin{array}{lll}
C^j(\Delta,\Lambda^p)_{\Q} & & i=0 \\
0                                 & & i\neq 0, \\
\end{array}\right.
\]
since
\[
H^{r-i}(\Star_{\sigma}(\Delta),\Q)=
\left\{\begin{array}{lll}
\Q & & i=0 \\
0  & & i\neq 0 \\
\end{array}\right.
\]
for all $\sigma\in\Delta(j)$.
Hence we get
\[
H_{\II}^j(H_{\I}^i(K^{\tcdot,\tcdot}_{\Q}))=
\left\{\begin{array}{lll}
H^j(\Delta,\Lambda^p)_{\Q} & & i=0 \\
0                                & & i\neq 0. \\
\end{array}\right.
\]
Consequently, by one of the two spectral sequences for the double complex,
we see that the associated single complex $(K^{\tcdot}_{\Q},d'+d'')$
has the cohomology group
\[
H^k(K^{\tcdot}_{\Q})=H^k(\Delta,\Lambda^p)_{\Q}
\qquad\mbox{for all}\quad k.
\]

On the other hand, for each fixed $i$, we have an isomorphism of complexes
\[
K^{i,\tcdot}_{\Q}=\bigoplus_{\varphi\in\Delta(r-i)}
\left(C^{\tcdot}(\Gamma_{\varphi},\Lambda^p)\otimes_{\Z}
(\det\varphi)^{-1}\right)_{\Q},
\]
the coboundary map for the left hand side being $d''$, while that for
the right hand side is $\delta\otimes\id$.
Thus by Proposition \ref{prop_vanishing_simplicialcone}, we get
\[
H_{\II}^j(K^{i,\tcdot}_{\Q})=
\left\{\begin{array}{lll}
\bigoplus_{\varphi\in\Delta(r-i)}
  \left((\bigwedge^p(M\cap\varphi^{\perp}))
  \otimes_{\Z}(\det\varphi)^{-1}\right)_{\Q} & & j=0 \\
0 & & j\neq 0, \\
\end{array}\right.
\]
hence
\[
H_{\I}^i(H_{\II}^j(K^{\tcdot,\tcdot}_{\Q}))=0
\qquad\mbox{for}\quad j\neq 0.
\]
Consequently, we have
\[
H^k(\Delta,\Lambda^p)_{\Q}=
H^k(K^{\tcdot}_{\Q})=H_{\I}^k(H_{\II}^0(K^{\tcdot,\tcdot}_{\Q}))
\]
for all $k$.
The extreme left hand side is nonzero only when $0\leq k\leq p$,
while the extreme right hand side is nonzero only when $p\leq k$.

The asserted perfect pairing is a consequence of the canonical
identification, as in \cite[Thm.\ 4.1]{lefsch}, of the ${\Q}$-dual of
$H^p(\Delta,\Lambda^p)_{\Q}=H_{\I}^p(H_{\II}^0(K^{\tcdot,\tcdot}_{\Q}))$
with $(\bigwedge\nolimits^r N_{\Q})\otimes_{\Q}
H^{r-p}(\Delta,\Lambda^{r-p})_{\Q}$.
The rest of the assertion follows from Theorem \ref{thm_algderham}.
\qed
\bigskip

The following is an important generalization, due to Ishida, of our earlier
result stated as Corollary \ref{cor_vanishing_convexsupport} later.

\begin{Theorem} \label{thm_vanishing_ishida}
{\rm (Ishida)}
Let $\Delta$ be a finite simplicial fan for $N\cong\Z^r$ which may not be
complete.
If there exist a finite complete simplicial fan $\tDelta$ and
a $\rho\in\tDelta(1)$ such that
$\Delta=\tDelta\setminus\Star_{\rho}(\tDelta)$, then
\[
H^q(\Delta,\Lambda^p)_{\Q}=0\qquad\mbox{for all}\quad q\neq p.
\]
Consequently, the corresponding toric variety $X:=T_N\emb(\Delta)$
satisfies $H^l(X\an,\C)=0$ for odd $l$, while
\[
H^{2p}(X\an,\C)=H^p(\Delta,\Lambda^p)_{\C}\qquad\mbox{for}\quad
0\leq p\leq r.
\]
\end{Theorem}

\Proof
Since $\Star_{\rho}(\tDelta)$ (resp.\ $\Delta$) is a star closed subset
(resp.\ a subcomplex) of $\tDelta$, we have an exact sequence of
complexes
\[
0\longrightarrow C^{\tcdot}(\Star_{\rho}(\tDelta),\Lambda^p)
\longrightarrow C^{\tcdot}(\tDelta,\Lambda^p)\longrightarrow
C^{\tcdot}(\Delta,\Lambda^p)\longrightarrow 0.
\]

Consider the projection $N\rightarrow\oN:=N/\Z(N\cap\rho)$.
For each $\sigma\in\Star_{\rho}(\tDelta)$,
its image $\osigma:=(\sigma+(-\rho))/\R\rho$ under the projection
$N_{\R}\rightarrow\oN_{\R}$ is a strongly convex rational polyhedral cone,
and
\[
\oSigma:=\{\osigma\;\mid\;\sigma\in\Star_{\rho}(\tDelta)\}
\]
is a finite complete simplicial fan for $\oN$.
The dual of $\oN$ is $\oM=M\cap\rho^{\perp}$, hence
$M\cap\sigma^{\perp}=\oM\cap\osigma^{\perp}$ holds for each
$\sigma\in\Star_{\rho}(\tDelta)$. Consequently, we have an isomorphism
of complexes
\[
C^{\tcdot-1}(\oSigma,\Lambda^{p-1})\stackrel{\sim}{\longrightarrow}
C^{\tcdot}(\Star_{\rho}(\tDelta),\Lambda^p).
\]
In view of Proposition \ref{prop_vanishing_complete} applied to
$\oSigma$, we thus have $H^q(\Star_{\rho}(\tDelta),\Lambda^p)_{\Q}=0$ for
$q\neq p$. Again by Proposition \ref{prop_vanishing_complete} applied
this time to $\tDelta$, we have $H^q(\tDelta,\Lambda^p)_{\Q}=0$ for $q\neq p$.
Hence we get $H^q(\Delta,\Lambda^p)_{\Q}=0$ for $q\neq p-1, p$ as well
as a long exact sequence
\[
0\longrightarrow H^{p-1}(\Delta,\Lambda^p)_{\Q}\longrightarrow
H^p(\Star_{\rho}(\tDelta),\Lambda^p)_{\Q}\longrightarrow
H^p(\tDelta,\Lambda^p)_{\Q}\longrightarrow H^p(\Delta,\Lambda^p)_{\Q}
\longrightarrow 0.
\]
In particular, the vanishing of $H^{p-1}(\Delta,\Lambda^p)_{\Q}$ is
equivalent to the injectivity of
\[
H^p(\Star_{\rho}(\tDelta),\Lambda^p)_{\Q}\rightarrow
H^p(\tDelta,\Lambda^p)_{\Q}.
\]

To show the latter, we now introduce another finite complete simplicial fan
$\tPhi$ for $N$ containing $\rho$ such that
$\Star_{\rho}(\tPhi)=\Star_{\rho}(\tDelta)$.
This $\tPhi$ and $\Phi:=\tPhi\setminus\Star_{\rho}(\tPhi)$ turn out to be
easier to handle, and we will be able to show $H^q(\Phi,\Lambda^p)_{\Q}=0$
for $q\neq p$ in Lemma \ref{lem_vanishing_P1bundle} below.
Hence by the same argument as above applied to $\tPhi$ instead of
$\tDelta$ we have the injectivity of
\[
H^p(\Star_{\rho}(\tPhi),\Lambda^p)_{\Q}\longrightarrow
H^p(\tPhi,\Lambda^p)_{\Q}.
\]
Clearly, there exists a finite simplicial complete fan $\tDelta'$ for
$N$ which is a subdivision of both $\tDelta$ and $\tPhi$
such that $\Star_{\rho}(\tDelta')=\Star_{\rho}(\tDelta)=\Star_{\rho}(\tPhi)$.

We claim the injectivity of the canonical homomorphisms
\[
H^p(\tDelta,\Lambda^p)_{\Q}\longrightarrow H^p(\tDelta',\Lambda^p)_{\Q}
\qquad\mbox{and}\qquad
H^p(\tPhi,\Lambda^p)_{\Q}\longrightarrow H^p(\tDelta',\Lambda^p)_{\Q}
\]
induced by the subdivisions. Consequently, we would get the injectivity of
\[
H^p(\Star_{\rho}(\tDelta),\Lambda^p)_{\Q}\rightarrow
H^p(\tDelta,\Lambda^p)_{\Q},
\]
since the canonical homomorphism from
\[
H^p(\Star_{\rho}(\tPhi),\Lambda^p)_{\Q}
=H^p(\Star_{\rho}(\tDelta),\Lambda^p)_{\Q}
=H^p(\Star_{\rho}(\tDelta'),\Lambda^p)_{\Q}
\]
to $H^p(\tPhi,\Lambda^p)_{\Q}$ is injective as above.

As for the proof of the above claim, it obviously suffices to prove the
injectivity of the scalar extension to $\C$ of the canonical homomorphisms
in question.
Let $\tX$, $\tX'$, $\tY$ be the complete toric varieties associated to the
fans $\tDelta$, $\tDelta'$, $\tPhi$, respectively, and denote by
$f:\tX'\rightarrow\tX$ and $g:\tX'\rightarrow\tY$ the equivariant proper
birational morphisms associated to the subdivisions of the fans.
By Theorem \ref{thm_algderham} and
Proposition \ref{prop_vanishing_complete}, we have
\begin{eqnarray*}
H^p(\tDelta,\Lambda^p)_{\C} &=& H^{2p}((\tX)\an,\C), \\
H^p(\tDelta',\Lambda^p)_{\C} &=& H^{2p}((\tX')\an,\C),\\
H^p(\tPhi,\Lambda^p)_{\C} &=& H^{2p}((\tY)\an,\C).
\end{eqnarray*}
The canonical homomorphisms in question coincide with
\[
f^*:H^{2p}((\tX)\an,\C)\rightarrow H^{2p}((\tX')\an,\C)
\qquad\mbox{and}\qquad
g^*:H^{2p}((\tY)\an,\C)\rightarrow H^{2p}((\tX')\an,\C),
\]
which are well-known to be injective by
$f_*f^*=\id$ and $g_*g^*=\id$, where $f_*$ and $g_*$ are the direct images
\[
f_*:H^{2p}((\tX')\an,\C)\rightarrow H^{2p}((\tX)\an,\C)
\qquad\mbox{and}\qquad
g_*:H^{2p}((\tX')\an,\C)\rightarrow H^{2p}((\tY)\an,\C).
\]

Here is how we define the new complete simplicial fan $\tPhi$ for $N$
which satisfies the required properties mentioned above:

$\Delta^{\flat}:=\{\tau\in\Delta\mid \tau+\rho\in\tDelta\}$ is easily
seen to be a fan for $N$. It can be thought of as the
``boundary'' of $\Delta$ as well as the ``link'' of $\rho$ in $\tDelta$
so that $\Star_{\rho}(\tDelta)=\{\tau+\rho\mid \tau\in\Delta^{\flat}\}$.
The projection $N\rightarrow\oN$ induces a bijection from each
$\tau\in\Delta^{\flat}$ to its image $\otau\in\oSigma$.
The cone $-\rho$ in $N_{\R}$ is not contained
in $\Star_{\rho}(\tDelta)$ nor $\Delta^{\flat}$. Hence
\[
\tPhi:=\{\tau+\rho\;\mid\;\tau\in\Delta^{\flat}\}\;{\textstyle \coprod}\;
\Delta^{\flat}\;{\textstyle \coprod}\;
\{\tau+(-\rho)\;\mid\;\tau\in\Delta^{\flat}\}
\]
is a finite simplicial complete fan for $N$ satisfying
$\Star_{\rho}(\tPhi)=\Star_{\rho}(\tDelta)$.
Let $\Phi^{\flat}:=\Delta^{\flat}$ and
\[
\Phi:=\Phi^{\flat}\;{\textstyle \coprod}\;
\{\tau+(-\rho)\;\mid\;\tau\in\Phi^{\flat}\}.
\]
We clearly have $\Phi=\tPhi\setminus\Star_{\rho}(\tPhi)$ and
$\Phi^{\flat}=\{\tau\in\Phi\mid\tau+\rho\in\tPhi\}$.

We are done in view of Lemma \ref{lem_vanishing_P1bundle} below
and Theorem \ref{thm_algderham}
\qed

\begin{Lemma} \label{lem_vanishing_P1bundle}
In the notation above, we have
$H^{\tcdot}(\Phi,\Lambda^p)_{\Q}=H^{\tcdot}(\oSigma,\Lambda^p)_{\Q}$.
In particular,
\[
H^q(\Phi,\Lambda^p)_{\Q}=0\qquad\mbox{for all}\quad q\neq p,
\]
hence the canonical homomorphism
$H^p(\Star_{\rho}(\tPhi),\Lambda^p)_{\Q}\rightarrow H^p(\tPhi,\Lambda^p)_{\Q}$
is injective.
\end{Lemma}

\Proof
As we mentioned above,
$\Phi^{\flat}=\Delta^{\flat}$ is in one-to-one correspondence with
$\oSigma$ by the map which sends $\tau\in\Phi^{\flat}$ to its isomorphic
image under the projection $N_{\R}\rightarrow\oN_{\R}=N_{\R}/\R\rho$.
On the other hand, $\Phi^{\flat}(q)$ for each $q$ is in one-to-one
correspondence with $\Star_{-\rho}(\Phi)(q+1)$ by the map which sends
$\tau\in\Phi^{\flat}$ to $\tau+(-\rho)\in\Star_{-\rho}(\Phi)$.
We have $\oM\cap\otau^{\perp}=M\cap(\tau+(-\rho))^{\perp}$ and an exact
sequence
\[
0\longrightarrow(\oM\cap\otau^{\perp})_{\Q}\longrightarrow
(M\cap\tau^{\perp})_{\Q}\longrightarrow\Q\longrightarrow 0,
\]
where the second arrow from the right is the interior product with the
unique primitive element $-n(\rho)$ of $N$ contained in the cone $-\rho$.
As a result, we have an exact sequence
\[
0\longrightarrow\bigwedge\nolimits^{p-q}(\oM\cap\otau^{\perp})_{\Q}
\longrightarrow\bigwedge\nolimits^{p-q}(M\cap\tau^{\perp})_{\Q}
\longrightarrow\bigwedge\nolimits^{p-q-1}(\oM\cap\otau^{\perp})_{\Q}
\longrightarrow 0
\]
for each $q$, hence an exact sequence of complexes
\[
0\longrightarrow C^{\tcdot}(\oSigma,\Lambda^p)_{\Q}\longrightarrow
C^{\tcdot}(\Phi^{\flat},\Lambda^p)_{\Q}\stackrel{\iota}{\longrightarrow}
C^{\tcdot}(\oSigma,\Lambda^{p-1})_{\Q}\longrightarrow 0.
\]
On the other hand, we have
\[
C^q(\Phi,\Lambda^p)=C^q(\Phi^{\flat},\Lambda^p)
\oplus C^q(\Star_{-\rho}(\Phi),\Lambda^p)
=C^q(\Phi^{\flat},\Lambda^p)\oplus C^{q-1}(\oSigma,\Lambda^{p-1})
\]
for each $q$. We see easily that
$C^{\tcdot}(\Phi,\Lambda^p)_{\Q}$ coincides with
the mapping cone of the surjective homomorphism
$\iota:C^{\tcdot}(\Phi^{\flat},\Lambda^p)_{\Q}\rightarrow
C^{\tcdot}(\oSigma,\Lambda^{p-1})_{\Q}$, whose kernel is
$C^{\tcdot}(\oSigma,\Lambda^p)_{\Q}$ as we saw above. Consequently,
we have
\[
H^{\tcdot}(\Phi,\Lambda^p)_{\Q}=H^{\tcdot}(\oSigma,\Lambda^p)_{\Q}.
\]
In particular, we have
$H^q(\Phi,\Lambda^p)_{\Q}=H^q(\oSigma,\Lambda^p)_{\Q}=0$ for $q\neq p$
by Proposition \ref{prop_vanishing_complete}.
\qed

\begin{Corollary} \label{cor_vanishing_convexsupport}
Let $\Delta$ be a finite simplicial fan for $N\cong\Z^r$ such that
its support $|\Delta|$ is convex of dimension $r$. Then
\[
H^q(\Delta,\Lambda^p)_{\Q}=0\qquad\mbox{for all}\quad q\neq p.
\]
Consequently, the corresponding toric variety $X:=T_N\emb(\Delta)$
satisfies $H^l(X\an,\C)=0$ for odd $l$, while
\[
H^{2p}(X\an,\C)=H^p(\Delta,\Lambda^p)_{\C}\qquad\mbox{for}\quad
0\leq p\leq r.
\]
\end{Corollary}

\Proof
By Proposition \ref{prop_vanishing_complete}, we may assume that
$\Delta$ is not complete. Then there exists
a primitive element $n^{\circ}\in N$
such that $-n^{\circ}$ is contained in the interior of $|\Delta|$.
Let $\rho:=\R_{\geq 0}n^{\circ}$.
Denote by $\Delta^{\flat}$ the subcomplex of $\Delta$ consisting of
those $\sigma\in\Delta$ which are contained in the boundary of the
convex cone $|\Delta|$. Obviously,
\[
\tDelta:=\Delta\;{\textstyle \coprod}\;
\{\sigma+\rho\;\mid\;\sigma\in\Delta^{\flat}\}
\]
is a finite complete simplicial fan for $N$ such that
$\Delta=\tDelta\setminus\Star_{\rho}(\tDelta)$. We are done by
Theorem \ref{thm_vanishing_ishida} and Theorem \ref{thm_algderham}.
\qed
\bigskip

$\tPhi$, $\Phi$ and $\Phi^{\flat}$ appearing in the proof of
Theorem \ref{thm_vanishing_ishida}, Lemma \ref{lem_vanishing_P1bundle}
and Corollary \ref{cor_vanishing_convexsupport} are of independent
interest. Namely, let $n_0$ be the primitive element in $N$
such that $\rho=\R_{\geq 0}(-n_0)$, and choose a splitting
$N\cong \oN\oplus\Z n_0$. Then there exists a function
$\eta:\oN_{\R}\rightarrow\R$, which is $\Z$-valued on $\oN$
and piecewise linear with respect to the complete nonsingular fan
$\oSigma$ for $\oN$ so that, in terms of the graph
$g:\oN_{\R}\rightarrow N_{\R}$ of $\eta$ defined by
$g(\on):=\on+\eta(\on)n_0$ for each $\on\in\oN_{\R}$, we have
\begin{eqnarray*}
\Phi^{\flat} &=& \{g(\osigma)\;\mid\;\osigma\in\oSigma\} \\
\Phi &=& \Phi^{\flat}
      \coprod\{\tau+\R_{\geq 0}n_0\;\mid\;\tau\in\Phi^{\flat}\} \\
\tPhi &=& \Phi
      \coprod\{\tau+\R_{\geq 0}(-n_0)\;\mid\;\tau\in\Phi^{\flat}\}. \\
\end{eqnarray*}
It is easy to see that
$T_N\emb(\tPhi)\rightarrow T_{\oN}\emb(\oSigma)$ is a $\bP_1(\C)$-bundle,
while $T_N\emb(\Phi)\rightarrow T_{\oN}\emb(\oSigma)$ and
$T_N\emb(\Phi^{\flat})\rightarrow T_{\oN}\emb(\oSigma)$ are
the associated $\C$-bundle and the associated $\C^*$-bundle, respectively.

As in \cite[Thm.\ 4.3 and Prop.\ 4.4]{lefsch} in different notation
and in Park \cite{park}, \cite{park2}, $\eta$ determines an element
$\oeta\in H^1(\oSigma,\Lambda^1)$, the multiplication by which induces
a homomorphism
\[
H^{p-1}(\oSigma,\Lambda^{p-1})_{\Q}\longrightarrow
H^p(\oSigma,\Lambda^p)_{\Q}
\]
for each $p$.
Then we have the following:

\begin{Corollary} \label{cor_vanishing_C*bundle}
In the notation above, we have
\[
H^q(\Phi^{\flat},\Lambda^p)=0\qquad\mbox{for}\quad q\neq p-1,p.
\]
Moreover, $H^{p-1}(\Phi^{\flat},\Lambda^p)_{\Q}$
(resp.\ $H^p(\Phi^{\flat},\Lambda^p)_{\Q}$) coincides with the
kernel (resp.\ cokernel) of the homomorphism
\[
H^{p-1}(\oSigma,\Lambda^{p-1})_{\Q}\longrightarrow
H^p(\oSigma,\Lambda^p)_{\Q}
\]
induced by multiplication of the element $\oeta\in H^1(\oSigma,\Lambda^1)$.
\end{Corollary}

%%%%%%%%%%%%%%%%%%%%%%%%%%%%%%%%%%%%%%%%%%%%%%%%%%%

\bigskip

\begin{flushleft}
\begin{sc}
Mathematical Institute \\
Faculty of Science \\
Tohoku University \\
Sendai 980 \\
Japan\bigskip \\
\end{sc}
e-mail: odatadao@jpntuvm0.bitnet
\end{flushleft}

\end{document}